\documentclass[reqno,11pt]{amsart}

\numberwithin{equation}{section}

\usepackage{verbatim} 
\usepackage{xspace} 
\usepackage{amssymb}
\usepackage{amsmath}
\usepackage{amsthm} %
\usepackage[latin1]{inputenc}
\usepackage{graphicx}
\usepackage{amscd,amssymb,euscript,graphics}
\makeindex

\newcommand{\Z}{\ensuremath{\mathbb{Z}}}

\newcommand{\F}{\ensuremath{\mathbb{F}}}
\newcommand{\N}{\ensuremath{\mathbb{N}}}

\newtheorem{theorem}{Theorem}[section]

\newtheorem{cor}[theorem]{Corollary}
\newtheorem{df}[theorem]{Definition}

\newcount\theTime
\newcount\theHour
\newcount\theMinute
\newcount\theMinuteTens
\newcount\theScratch
\theTime=\number\time
\theHour=\theTime
\divide\theHour by 60
\theScratch=\theHour
\multiply\theScratch by 60
\theMinute=\theTime
\advance\theMinute by -\theScratch
\theMinuteTens=\theMinute
\divide\theMinuteTens by 10
\theScratch=\theMinuteTens
\multiply\theScratch by 10
\advance\theMinute by -\theScratch

\def\today{{\number\day\space
 \ifcase\month\or
  January\or February\or March\or April\or May\or June\or
  July\or August\or September\or October\or November\or December\fi
 \space\number\year}}

\begin{document}

\title[Lack of Hyperbolicity in Erd\"os--Renyi Random Graphs]{Lack of Hyperbolicity in Asymptotic Erd\"os--Renyi Sparse Random Graphs}

\author{Onuttom Narayan, Iraj Saniee and Gabriel H. Tucci} 
\thanks{I. Saniee and G. H. Tucci are with Bell Laboratories, Alcatel-Lucent, 600 Mountain Avenue, Murray Hill, New Jersey 07974, USA. 
O. Narayan is at Department of Physics, University of California, Santa Cruz, California 95064, USA.
E-mail: narayan@wagner.ucsc.edu, iis@research.bell-labs.com, gabriel.tucci@alcatel-lucent.com}

\begin{abstract}
In this work we prove that the giant component of the Erd\"os--Renyi random graph $G(n,c/n)$ for $c$ a constant greater than $1$ (sparse regime), is not Gromov $\delta$--hyperbolic for any $\delta$ with probability tending to one as $n\to\infty$. As a corollary we provide an alternative proof that the giant component of $G(n,c/n)$ when $c>1$ has zero spectral gap almost surely as $n\to\infty$.
\end{abstract}

\maketitle

\section{Introduction and Motivation}

\noindent Random graphs constitute an important and active research area
with numerous applications to geometry, percolation theory, information
theory, queuing systems and communication networks, to mention a
few.  They also provide analytical means to settle prototypical
questions and conjectures that may be harder to resolve in specific
circumstances (such as statistical evidence for hyperbolicity or its lack 
via curvature plots, as discussed in \cite{NS2}, which is our focus here).  
In this work we study two questions regarding the
asymptotic geometry of Erd\"os--Renyi random
graphs ~\cite{Erdos1,Gilbert,Chung}, partly motivated by inference
that random graphs may be hyperbolic \cite{jonck} or may have 
a spectral gap \cite{coja}.   These and other authors use the term random graph
in different senses. To fix definition and notation, we call 
$G(n,p_{n})$ a random graph where $n$ is the the number of nodes and
$p_{n}$ is the probability of an edge between any node pair, independently of all other
edges.  The construction of a $G(n,p_{n})$ consists of connecting any pair of 
these $n$ nodes independently with probability $p_{n}$\footnote{This is
actually the $G(n,p)$ model of a random graph due to
Gilbert ~\cite{Gilbert}, rather than the Erd\"os--Renyi ~\cite{Erdos1}
model known as $G(n,M)$; but we follow the now almost universal
proclivity of referring to these as Erd\"os--Renyi random graphs.}.

\vspace{0.3cm}
\noindent Our main result is that in the constant average--degree
regime $p_{n}=c/n$ with $c$ a constant greater than $1$, with probability 
approaching one these graphs
are not $\delta$--{\em hyperbolic} in the sense of Gromov ~\cite{Gromov} (which we make precise in Section~\ref{non-hyper}) 
for any non--negative $\delta$. One might think that this is equivalent to the lack of spectral gap,
since Gromov's notion of hyperbolicity and the linear isoperimetric
inequality are intimately related in a coarse sense, see
~\cite{bridson}.  In fact, despite the connection between the
two, neither one implies the other as we discuss in
more detail in Section~\ref{non-hyper}. This implies that the questions
of hyperbolicity and spectral gap of random graphs need to be
addressed independently.

\vspace{0.3cm}
\noindent This paper is organized as follows. In Section \ref{non-hyper}, we show that the giant component of $G(n,c/n)$ is not $\delta$--hyperbolic for any $\delta$ with probability tending to one as $n\to\infty$. This implies that for every positive $\delta$ there are triangles in $G(n,c/n)$ that are not $\delta$--thin. These triangles are called $\delta$--fat.  In Section \ref{sim}, we present plots that suggest that ``fat'' triangles not only exist almost surely as $n\to\infty$ but are abundant in these random graphs. We also present numerical results that show a surprising degree of closeness between the spectral distribution of the normalized Laplacian of $G(n,c/n)$ and that of $c$--regular trees using well--known explicit formulas due to McKay ~\cite{McKay}.

\section{Non--hyperbolicity for the Erd\"os--Renyi Random Graphs}
\label{non-hyper}
\subsection{Relationship Between Hyperbolicity and Spectral Gap}
\label{sec:nonequiv}
It is known that the giant component of $G(n,c/n)$ does not have a spectral gap 
almost surely in the regime $c>1$.  This follows, for instance, from \cite{Itai} to cite a
recent paper. This means that as 
$n\to\infty$, the smallest non--zero eigenvalue of the Laplacian of the giant component of $G(n,c/n)$ (see \cite{Chung}) goes to zero.
One natural way to prove this is to show that there are arbitrarily long paths 
with a unique single attachment to the graph almost surely. In this Section, we prove a stronger result
for the Erd\"os--Renyi random
graphs in the sparse regime ($p_n=c/n,c>1$), showing that these graphs possess arbitrarily 
long loops with the ends attached to the rest of the graph, thus demonstrating that these graphs 
are not $\delta$--hyperbolic in the sense of Gromov ~\cite{Gromov} for any $\delta$.  

\vspace{0.3cm}
\noindent To be more precise about the expression ``$\delta$--hyperbolic in the sense of Gromov'' for a
family of finite graphs, let $G=(V,E)$ be a 
(finite) graph together with an edge metric $d$ (thus $d$ satisfies
the triangle inequality).  Let $[XY]$ denote a shortest path between vertices $X$ and $Y$ in $G$.
A triangle with vertices $X,Y$ and $Z$ is said to be $\delta$--thin if 
\begin{equation}
[XY] \subseteq \mathcal{N}([YZ],\delta) \cup \mathcal{N}([ZX],\delta) 
\end{equation}
where $\mathcal{N}([XY],\delta)$ is the $\delta$ neighbourhood of $[XY]$ and so on.
The graph $G$ is said to be $\delta$--hyperbolic if all its triangles are $\delta$--thin.
%
%
%
%
Intuitively, $\delta$--hyperbolicity means that any three shortest paths $[XY]$, $[YZ]$ and $[ZX]$
between any triple of vertices $X,Y$ and $Z$ in $G$ come to within a distance $\delta$ of each other for
a some fixed $\delta \geq 0$.  
Thus trees are $0$--hyperbolic, the two dimensional square grid is not $\delta$--hyperbolic for any finite 
$\delta$ but any finite graph with diameter $\Delta$ is $\Delta$--hyperbolic.

\vspace{0.3cm}
\noindent We say a family of graphs $\{G_n\,:\,n \geq 1 \}$ is $\delta$--hyperbolic if each member $G_n$ is $\delta$--hyperbolic
for a fixed value $\delta>0$. We say a family is asymptotically $\delta$--hyperbolic if for large enough $n$ all $G_n$ are $\delta$--hyperbolic.
When a family $\{G_n\,:\,n \geq 1 \}$ is not $\delta$--hyperbolic for any $\delta \geq 0$, then it must be
the case that for any positive $\delta$ there is an $n$ such that there are some $\delta$--fat triangles in $G_n$.
This is precisely in the sense in which we prove that the family $\{H_{n}\,:\,n\geq 1\}$ where $H_{n}$ is the giant component of $G(n,p_n = c/n)$ with $c>1$ is not $\delta$--hyperbolic.

\vspace{0.3cm}
\noindent The concept of hyperbolicity is usually associated with the existence of a 
spectral gap. This is because for standard hyperbolic spaces with constant negative curvature,
the eigenvalues of the Laplace operator are bounded away from zero \cite{Chavel}.  
\noindent Indeed, it might be thought that the existence of a spectral gap and 
$\delta$--hyperbolicity are equivalent: the first is clearly equivalent 
to the existence of a linear isoperimetric inequality, and the second is shown to 
be equivalent to a linear isoperimetric inequality, for example see Proposition III.2.7 
of ~\cite{bridson}. However, the term ``linear 
isoperimetric inequality" is used in different senses in the two cases. In 
the first, the entire perimeter of any arbitrary subset $S$ has to be 
considered. In the second, disk--like subsets are considered, and only the 
loop part of the perimeter (ignoring any boundary edges on the ``flat" part 
of the disk) is used. Thus neither does the existence of a spectral gap imply
$\delta$--hyperbolicity nor vice versa.\footnote{We thank M.R. Bridson for a 
useful discussion on this point.} 

\vspace{0.3cm}
\noindent As examples to illustrate this fact, we note that a graph that
consists of an infinite chain (the integers $\Z$) has a zero Cheeger
constant and a zero spectral gap, even though it is $\delta$--hyperbolic
because -- as discussed in the next subsection -- all tree graphs
trivially are. On the other hand, the Cayley graph associated with
the product of two free groups, $G=\F_{2}\times\F_{2},$ has a
positive Cheeger constant and non-zero spectral gap.  But since it
includes the graph $H = \Z\times\Z$ (the Euclidean grid) as a
subgraph it is not hyperbolic.  Thus questions regarding the spectral
gap and hyperbolicity need to be addressed independently.

\subsection{Positive Measure of Large Loops}
\begin{figure}[htp]
\centering
\includegraphics[width=6cm, height=6cm]{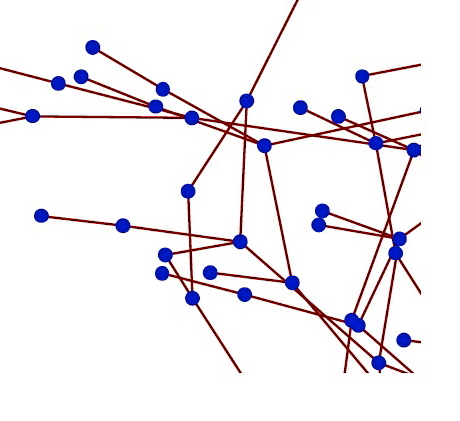}
\caption{A small segment of a random graph $G(n,2/n)$ viewed close up for $n=1000$.}
\label{fig:segment}
\end{figure}

It is commonly stated that the Erd\"os--Renyi random graphs
are ``tree--like'' for large values of $n$, on the strength of the
notion that any small neighborhood (the ``small scale'') has a very
small probability of localized links, see for example
figure~\ref{fig:segment}  (see ~\cite{MM09, karrer-newman}).
This ``treeness'' in the small scale is sometimes loosely interpreted
to imply that random graphs are {\em hyperbolic}. There are several
concerns about these heuristic notions and clarification is needed.
First, the probability regime of the construction of the random graph needs to be specified.  Second, more formal definitions
of small, middle and large scale are needed. As it is well known, there are different regimes in the $G(n,p)$ model of a random
graph:
\begin{enumerate}
\item $p = o(1/n)$, then the random graph is a large collection of disconnected small trees.

\item $p = c/n$ with $0<c<1$, then all the connected components of the graph are either trees or unicycle components. The giant connected component is a tree and has $O(\log(n))$ nodes.

\item $p = c/n$ with $c>1$, then a giant component emerges. This one has $\gamma(c)n$ nodes almost surely where $\gamma$ a function depending on $c$ and independent on $n$. Also the average degree of a node is bounded away from $0$.

\item $p = c\log(n)/n$ with $c>1$, then the graph is almost surely connected.  
\end{enumerate}

\vspace{0.3cm}
\noindent Beyond these, for example when $p=\Omega(1/\log(n))$, there is a single highly connected component whose average nodal degree
is unbounded as $n\to\infty$.

\begin{figure}[htp]
\centering
\includegraphics[width=12cm, height=9cm]{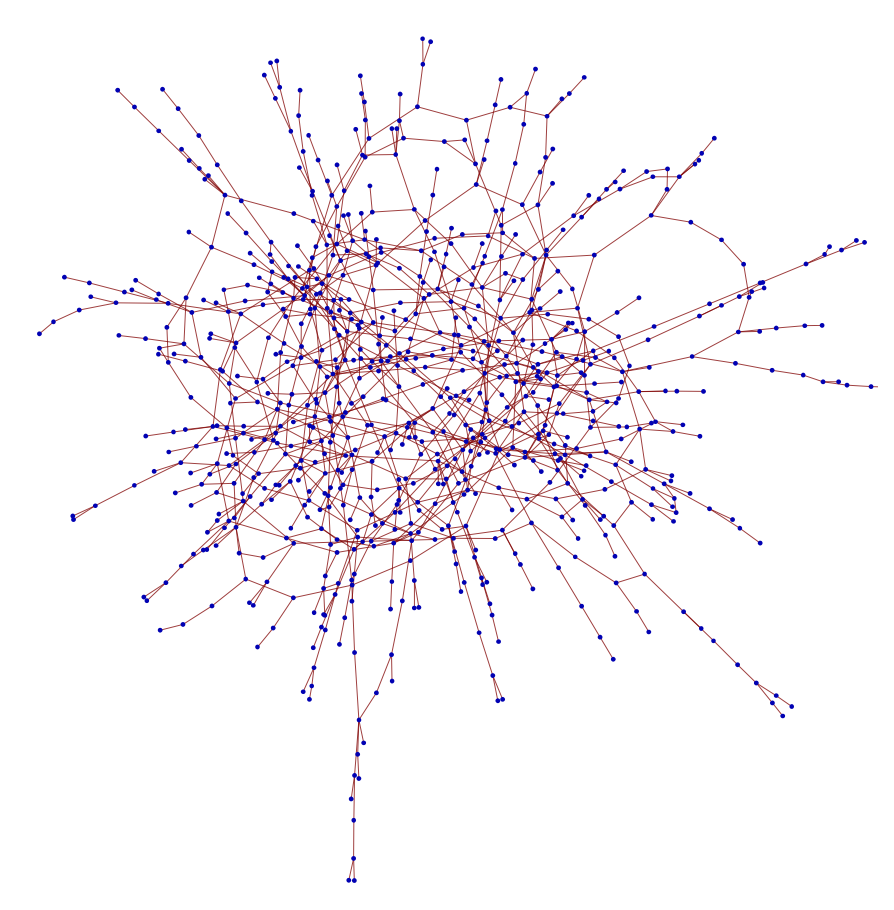}
\caption{A random graph $G(n,2/n)$ with $n=1000$. There are
loops of all sizes up to $O(\log(n))$, the order of its diameter.}
\label{fig:visual}
\end{figure}

\vspace{0.3cm}
\noindent With these clarifications, we make the following observations. First, random graphs in the $p=c/n$ (middle) regime are not
$\delta$--hyperbolic, in the sense that they contain $\delta$--fat triangles for arbitrary large $\delta$ almost surely as $n \rightarrow \infty$.  This is proved in Theorem~\ref{ourtheorem}. This observation was made experimentally in ~\cite{NS2} (see
the taxonomy chart there) and for which we provide a proof in this work.  Figure~\ref{fig:visual} provides a visualization of this
claim. Second, simulations suggest that the proportion of $\delta$--fat triangles is not only positive but is in fact quite significant for
any $\delta$ as $n$ grows. These are shown in Section~\ref{sec:curvplots}.

\begin{df}
Let $\{G_{n}\}_{n=1}^{\infty}$ be a family of random graphs. We say that a property holds asymptotically almost surely if the probability $p_{n}$ of this to occur goes to one as $n\to\infty$.
\end{df}

\begin{theorem}
\label{ourtheorem}
For every non--negative $\delta$, the giant component of $G(n,c/n)$ with $c>1$ is not $\delta$--hyperbolic asymptotically almost surely.
\end{theorem}

\begin{proof}

Let $G=G(n,c/n)$ be an Erd\"os--Renyi random graph with $c>1$ and let $H$ be its giant connected component. It is well known that $H$ has $\gamma(c)n$ nodes asymptotically almost surely where $\gamma(c)$ is a function on $c$ and independent on $n$ (see \cite{Chung} for a proof of this result). Take $\rho>0$ and let us expose $(1-\rho)n$ of the nodes of $G$ (by exposing we mean to generate the Erd\"os--Renyi random graph generated by these nodes). We call this set the {\it exposed} set. The remaining set of $\rho n$ is called the {\it hidden} set and is denoted by $\mathcal{H}$. It is easy to see that if $\rho$ is such that $0<\rho<\frac{c-1}{c}$ then the exposed set has a giant connected component of size $(1-\rho)\gamma\big((1-\rho)c\big)n$. Moreover, by taking $\rho$ as before we see that the giant component of the exposed set is contained in the giant component of the whole graph $G$. This is because in the graph $G$ there is a unique component of size proportional to $n$ all the other components have size $O(\log(n))$.

\vspace{0.3cm}
\noindent Let $v$ be a node in $\mathcal{H}$ and let $k$ be a positive integer. The probability that $v$ has only two neighbors in $\mathcal{H}$ and no other neighbor is equal to $(\rho n-1)(\rho n-2)p^{2}(1-p)^{n-3}/2$. This probability converges to $\rho^2 c^{2}e^{-c}/2$ as $n\to\infty$. The probability of their neighbors to have another unique neighbor in $\mathcal{H}$ is asymptotically $\rho ce^{-c}$. Moreover, the probability of the following neighbors to have a unique neighbor and so on until the nodes $k$ and $k'$ (see figure \ref{fig:handle}) are in the giant component of the exposed set is 
\begin{equation}\label{p_k}
p_{k} \approx \frac{\rho^{2k}(c e^{-c})^{2k+2} \gamma\big((1-\rho)c\big)^2(1-\rho)^2}{2e^{-c}}.
\end{equation}
Where the symbol $\approx$ denotes that the quantities are asymptotically equal. 

\vspace{0.3cm}
\noindent We say that $v$ is the base of a $k$--handle if $v$ is a node as in figure \ref{fig:handle}. Let $k$ and $k'$ be the nodes in the $k$--handle that belong  to the giant component of the exposed set. Note that since the nodes $k$ and $k'$ were already in the giant component of the exposed set there exists at least one shortest path connecting them inside the exposed set and not passing through the node $v$. Let $X_{v}$ be the random variable that is equal to $1$ if the node $v$ is the base of a $k$--handle with $v\in\mathcal{H}$ and $0$ otherwise. 
\begin{figure}[htp]
\centering
\includegraphics[height=6cm]{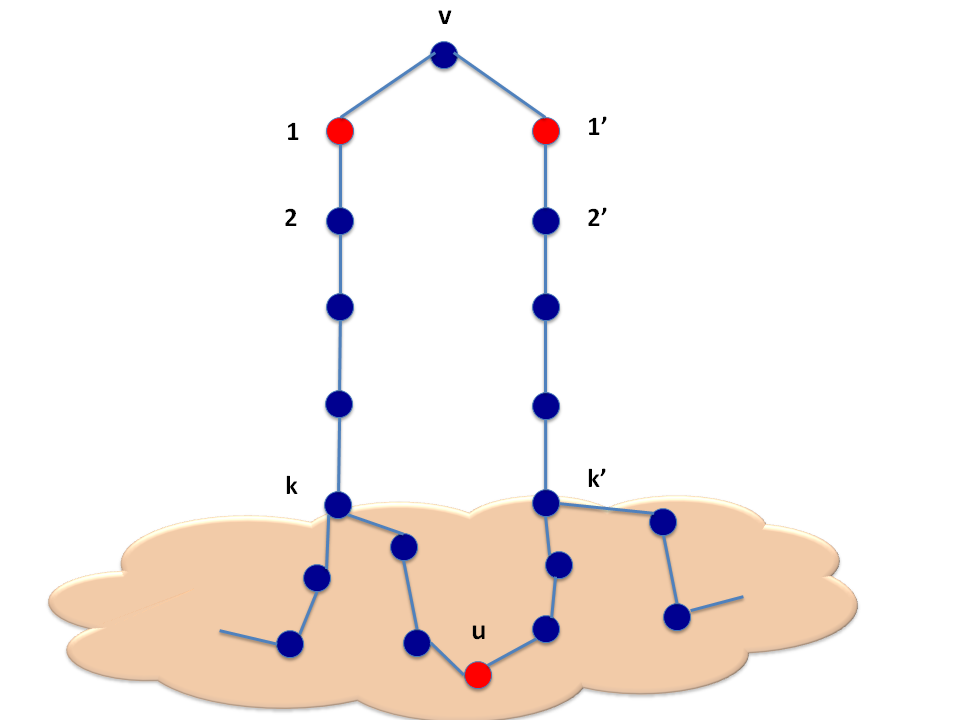}
\caption{Depiction of a $k$--handle based at the point $v$. The nodes $k$ and $k'$ are the nodes in the $k$--handle that belong to giant component of the exposed set. The node $u$ is the midpoint of $k$ and $k'$ for the shortest path connecting these two nodes without passing through the node $v$.}
\label{fig:handle}
\end{figure}
Let $1$ and $1'$ be the neighbors of $v$ and let $u$ be the midpoint of the points $k$ and $k'$ (these are the nodes marked in red in the figure) in any path that connects $k$ and $k'$ without passing through $v$. It is clear that the geodesic triangle $\Delta(11'u)$ is at least $\lfloor k/2\rfloor$--fat. 

\vspace{0.3cm}
\noindent Define $Y_{k}$ to be the random variable $Y_{k}=\sum_{v\in\mathcal{H}}{X_{v}}$. To prove the existence of a $\lfloor k/3\rfloor$--fat triangle in the giant component almost surely it is enough to prove that $\mathbb{P}(Y_{k}\geq 1)\to 1$ as $n\to\infty$. Moreover, we will prove the following stronger result. For every constant $0<t<(2(c-\log(\rho)))^{-1}$, there are almost surely $(t\log(n))$--handles as $n\to\infty$. Taking $k=t\log(n)$ in equation (\ref{p_k}) we obtain that  
\begin{eqnarray}
p_{t\log(n)} &\approx & \Bigg(\frac{\gamma((1-\rho)c)^2 c^2 (1-\rho)^2 e^{-c}}{2}\Bigg) (\rho c e^{-c})^{2t\log(n)},\\
&=& \Bigg(\frac{\gamma((1-\rho)c)^2 c^2 (1-\rho)^2 e^{-c}}{2}\Bigg) n^{2t\log(\rho ce^{-c})}.
\end{eqnarray}
Let us define 
$$
\theta(c):= \Bigg(\frac{\gamma((1-\rho)c)^2 c^2 (1-\rho)^2 e^{-c}}{2}\Bigg).
$$
Then the expected number of $(t\log(n))$--handles in $\mathcal{H}$ is
\begin{eqnarray}
\mathbb{E}(Y_{t\log(n)}) & = & \sum_{v\in\mathcal{H}}{\mathbb{E}(X_{v})} = \rho n p_{t\log(n)} \approx \rho\theta(c)n^{1+2t\log(\rho ce^{-c})}, \\
& = & \rho\theta(c)n^{1+2t\log(c)-2t\log(e^{c}/\rho)},
\end{eqnarray}
since the random variables $X_{v}$ are identically distributed and Bernoulli. Since by assumption $0<t<(2(c-\log(\rho)))^{-1}$, we see that this quantity goes to infinity as $n$ increases. Let $v$ and $w$ be two different nodes in $\mathcal{H}$ and let $q_{k}$ be the probability of having a $k$--handle based at $v$ and another based at $w$. It is rather easy to show that the quotient $q_{k}/p_{k}^{2}\to 1$ as $n\to\infty$. Recall that two random variables $X$ and $Y$ are independent if and only if $\mathbb{E}(X^{n}Y^{m})=\mathbb{E}(X^{n})\mathbb{E}(Y^{m})$ for all $n$ and $m$ integers greater or equal than $1$. Since $X_{v}^{n}=X_{v}$ for all $n\geq 1$ then to prove that $X_{v}$ and $X_{w}$ are asymptotically independent it is enough to show that $\mathbb{E}(X_{v}X_{w})/\mathbb{E}(X_{v})^2\to 1$. On the other hand,
\begin{equation}
\frac{\mathbb{E}(X_{v}X_{w})}{\mathbb{E}(X_{v})^2}=\frac{q_{k}}{p_{k}^2}\to 1.
\end{equation}
Therefore, we showed that $X_{v}$ and $X_{w}$ are asymptotically independent for all $v\neq w\in\mathcal{H}$. 

\vspace{0.3cm}
\noindent It is a straightforward calculation to show that the variance $\mathbb{V}(Y_{t\log(n)})$ satisfies 
\begin{equation}
\mathbb{V}(Y_{t\log(n)}) = \rho n p_{t\log(n)}(1- p_{t\log(n)}),
\end{equation}
since the random variables $X_{v}$ and $X_{w}$, as we showed, are asymptotically independent. Note that the probability of not having a $(t\log(n))$--handle based at $\mathcal{H}$ is equal to the probability of $Y_{t\log(n)}=0$. Hence, by Chebyshev's inequality 
\begin{eqnarray*}
0\leq \mathbb{P}(Y_{t\log(n)}=0) &\leq& \mathbb{P}\Bigg(|Y_{t\log(n)} - \mathbb{E}(Y_{t\log(n)})|\geq \frac{\mathbb{E}(Y_{t\log(n)})}{2}\Bigg)\leq \frac{4\mathbb{V}(Y_{t\log(n)})}{\mathbb{E}(Y_{t\log(n)})^{2}} \\
& = & \frac{4\rho n p_{t\log(n)}(1- p_{t\log(n)})}{\rho^2 n^2 p_{t\log(n)}^2} = \frac{4(1- p_{t\log(n)})}{\rho n p_{t\log(n)}}.
\end{eqnarray*}
By our election of $t$ we know that $p_{t\log(n)}\to 0$ and that $np_{t\log(n)}\to\infty$ as $n\to\infty$. Therefore, $\mathbb{P}(Y_{t\log(n)}=0)\to 0$ and our result follows.
\end{proof}

\noindent As a corollary of the previous proof we have an alternative proof of the following known result (e.g., see \cite{Itai}).

\begin{cor}
The giant component of $G(n,c/n)$ with $c>1$ has no spectral gap asymptotically almost surely.
\end{cor}

\noindent This result follows because the previously constructed $(t\log(n))$--handles are cut sets with $2t\log(n)+1$ nodes and only two boundary nodes.

\section{Simulations}\label{sim}
\subsection {Numerical Results on Percentage of Fat Triangles}
\label{sec:curvplots}
\begin{figure}[htp]
\centering
\includegraphics[width=7cm]{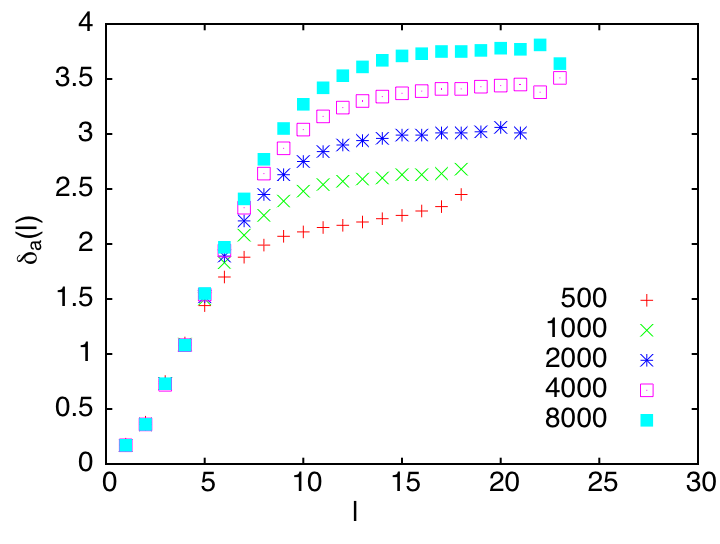}
\includegraphics[width=7cm]{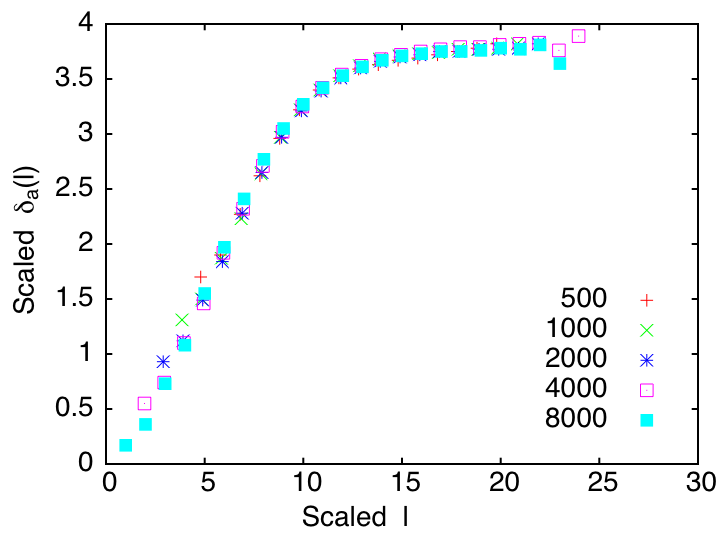}
\caption{Curvature plot for random graphs with $p=2/n$ and various
values of $n$. Only the giant component of each graph was retained,
and an average over many randomly chosen triangles in 40 instantiations
of the graph was performed. The right panel shows the same curves
as on the left, but shifted down and to the left by amounts
proportional to $\ln (n/8000).$}
\label{fig:curvplot}
\end{figure}

\noindent We have seen that random graphs $G(n,p_{n})$ in the regime $p_{n} = c/n$ are almost surely asymptotically non--hyperbolic. However, these random graphs appear to be non--hyperbolic in a much stronger sense. To see how, consider the chart in figure~\ref{fig:curvplot}.
This is an example of a curvature plot (see ~\cite{NS2}). For
any triangle $\Delta = ABC,$ the corresponding $\delta_\Delta$ is
defined by
\begin{equation}
\delta_\Delta = \min_D \max\Big\{ d(D;AB), d(D;BC), d(D;AC)\Big\}
\end{equation}
where $d(D;AB)$ is the distance between $D$ and the node on $AB$
that it is closest to. It can be shown (see ~\cite{bridson})
that the maximum of $\delta_\Delta$ over all possible triangles in
a graph is finite if and only if the graph is $\delta$--hyperbolic.
Instead of the maximum, figure~\ref{fig:curvplot} shows the average
value of $\delta_\Delta$ for all triangles whose shortest side is
$l,$ as a function of $l.$ Results for random graphs of various
sizes, with $p=2/n$ are shown. The results show a linear
increase in $\delta_a(l)$ saturating at a plateau whose height
increases as the size of the graph is increased. The same results
are rescaled in the right panel, where all the curves are shifted
down and to the left by amounts proportional to $\ln n.$ Thus the
plot shows $\delta_a(l) - c_1 \ln n$ versus $l - c_2 \ln n,$ where
$c_1$ and $c_2$ are adjusted to achieve the best possible fit. As
shown in the figure, except the leftmost part of each curve where
$l\sim O(1)$, all the curves for different $n$ collapse onto a
single universal curve. This, together with the fact that the curves
in the left panel also coincide before their plateaus implies that
the rising part of the universal curve is linear.  If $n\rightarrow\infty$,
any finite $l\gg1$ is on the rising part of the universal curve,
and therefore $\delta_a(l)$ increases linearly with $l,$ with the
plateau pushed out to $l\rightarrow\infty.$

\vspace{0.3cm}
\noindent Thus we see that a significant fraction of triangles in a typical
instantiation of $G(n,c/n)$ are $\delta$--fat, a stronger demonstration of non--hyperbolicity. Therefore, it seems that $\delta$--fat triangles not only exist almost surely but they are abundant! Even though we do not yet
have direct proof of this observation, figure~\ref{fig:curvplot}
clearly shows the predominance of fat triangles in $G(n,c/n)$
due to the increasing (average) $\delta$. Thus $G(n,c/n)$
random graphs are far from hyperbolic, contrary to folklore.

\begin{figure}[htp]
\centering
\includegraphics[width=8cm]{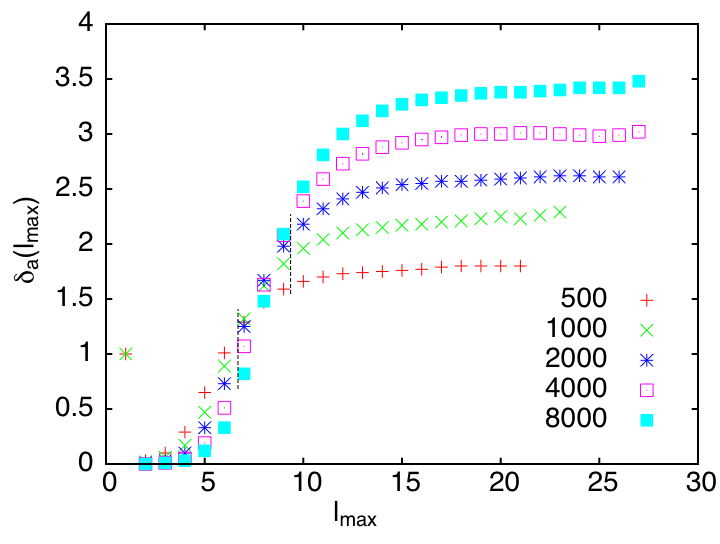}
\caption{Average $\delta_\Delta$ for all triangles with longest
side $l_{\max}$ as a function of $l_{\max}$, for various graph
sizes. The points in the region between the dashed lines show that,
as $n$ is increased, the range of $l_{\max}$ over which the curves
move downwards expands to the right.}
\label{fig:longside}
\end{figure}

\begin{figure}[htp]
\centering
\includegraphics[width=7cm]{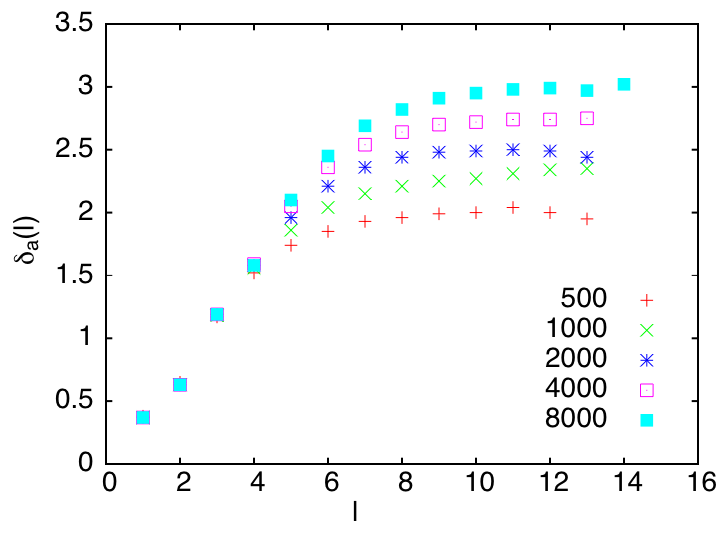}
\includegraphics[width=7cm]{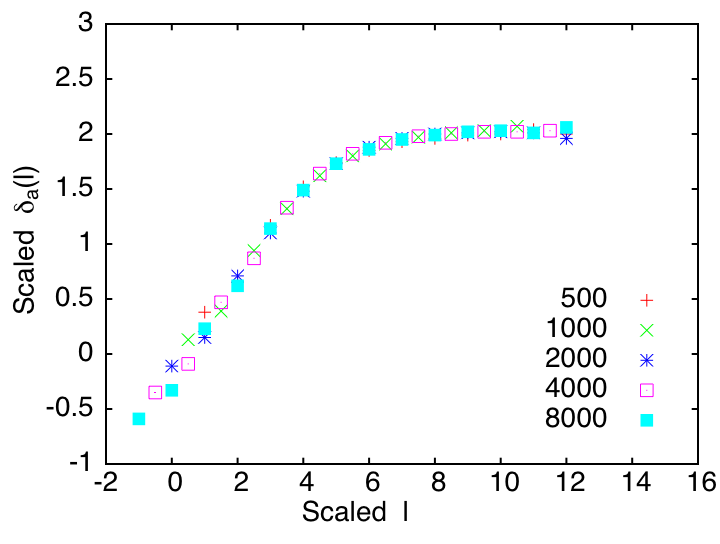}
\caption{Plots for random graphs with $p=3/n$, similar to
figure~\ref{fig:curvplot} for $p=2/n$.}
\label{fig:curvplot3}
\end{figure}

\vspace{0.3cm}
\noindent Figure~\ref{fig:longside} shows that if the average $\delta_\Delta$ for
all triangles with the same {\it longest\/} side $l_{\max}$ is plotted
as a function of $l_{\max},$ as $n$ increases, the height of the curves
decreases for small $l_{\max}$ and increases for large $l_{\max},$ with the
boundary between the two regions shifting to the right as $n$ increases.
Thus $\lim_{n\rightarrow\infty} \delta_a(l_{\max}) = 0$ for any fixed
$l_{\max},$ in accordance with the local tree-like structure. Similar 
results are seen for $p=3/n$ in figure~\ref{fig:curvplot3}.

\subsection{Some Simulation Results on the Bulk Region of the Spectrum}\label{bulk}

\begin{figure}[Htp]
\begin{center}
	\includegraphics[width=7cm, height=6cm]{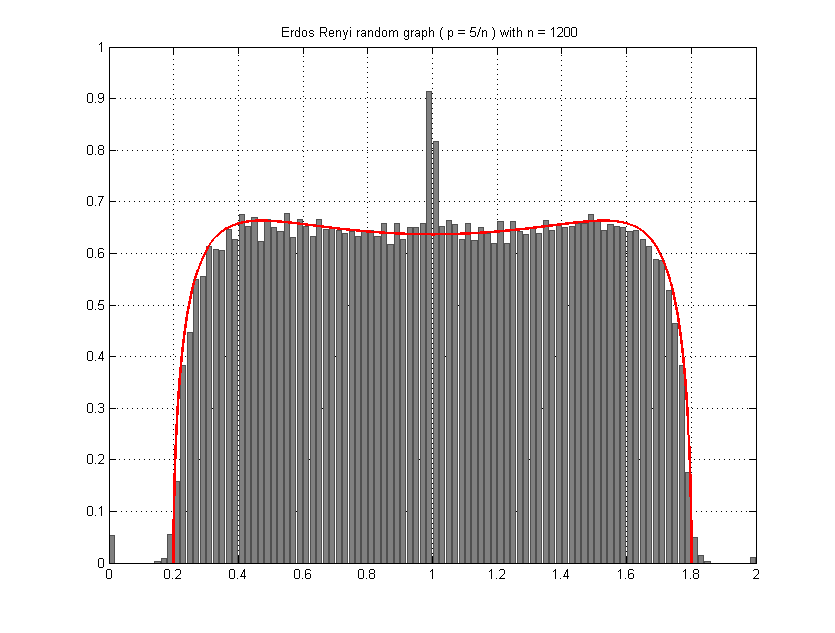}
	\includegraphics[width=7cm,height=6cm]{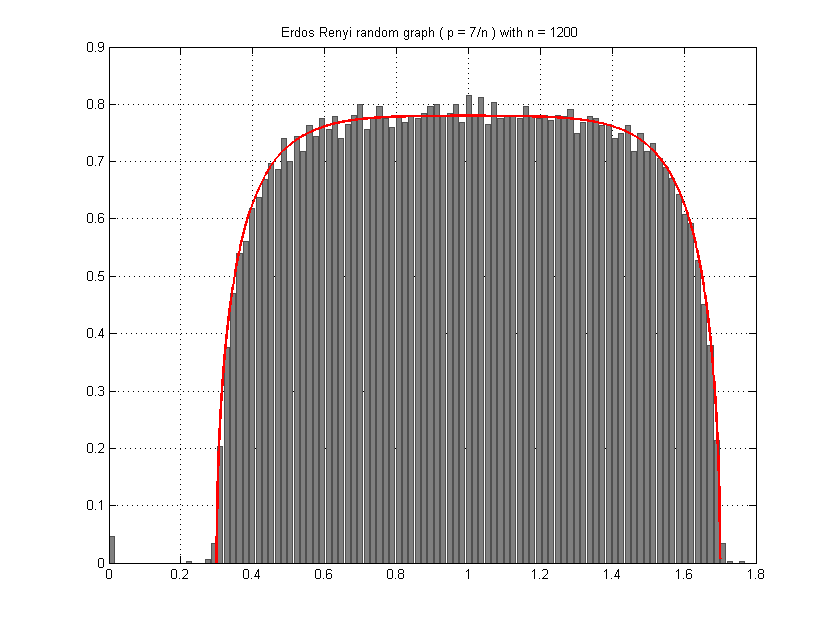}
\end{center}

\caption{The left chart in grey is $\mu_{n}$, the spectral density
of $G(n,5/n)$ for $n=1200$.  The right chart is $\mu_{n}$,
the spectral density of $G(n,7/n)$ for $n=1200$.  The red
curves are the McKay densities, the spectrum of the infinite regular
tree with degree $c=5$ and $c=7$ respectively.}
\label{fig:mckay1}
\end{figure}

\begin{figure}[Htp]
\begin{center}
	\includegraphics[width=7cm,height=6cm]{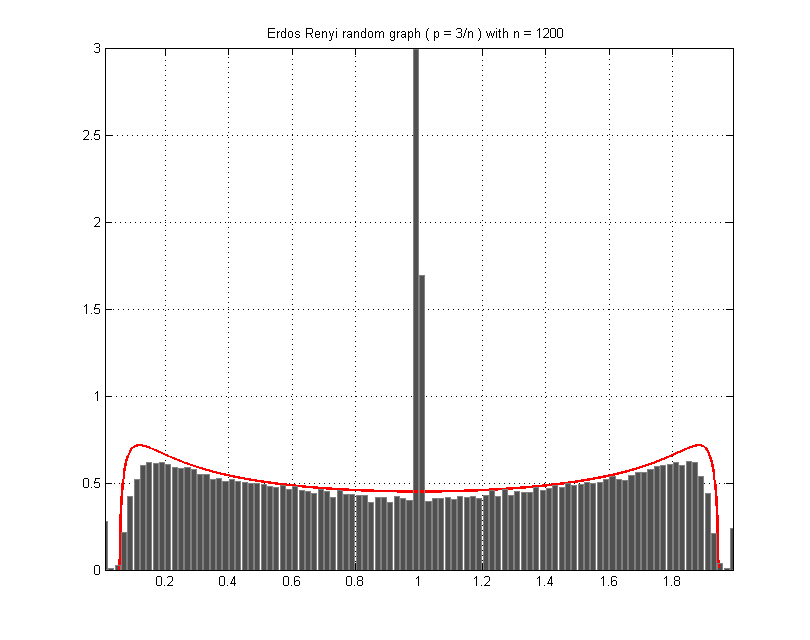}
	\includegraphics[width=7cm,height=6cm]{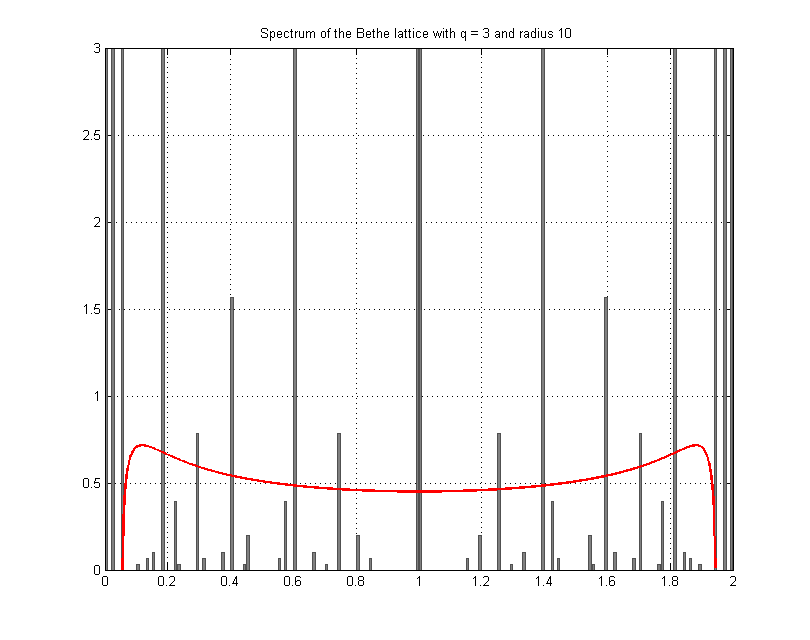}
\end{center}

\caption{The left chart in grey is $\mu_{n}$, the spectral density
of $G(n,3/n)$ for $n=1200$.  The red curve is the McKay
density, the spectral density of the infinite regular tree, with
degree $c=3$. In the right chart in grey is $\nu_{n}$, the spectral
density of the truncated tree with degree $c=3$ and radius $10$.
The red curve is the McKay density of degree $c=3$.}
\label{fig:mckay2}
\end{figure}

\noindent Here we present some simulations of the spectral measure $\mu_{n}$ for the Laplacian of the Erd\"os--Renyi graphs $G(n,c/n)$. It is known that these measures converge weakly to a probability measure $\mu_{\infty}$ (see \cite{Bord}). Observe in figures~\ref{fig:mckay1} and \ref{fig:mckay2} how close
these probabilities are in the bulk region to the McKay probability
measure, the spectral measure of the Laplacian of the $c$--regular
tree, that is given by ~\cite{McKay} to be
\begin{equation}\label{Mckay}
\mu(dx) = \frac{\sqrt{4(c-1)-c^2(1-x)^2}}{2\pi c(1-(1-x)^2)}
\cdot\mathbf{1}_{\big[1-\frac{2\sqrt{c-1}}{c},1+\frac{2\sqrt{c-1}}{c}\big]}.
\end{equation}

\vspace{0.3cm}
\noindent It is interesting to compare these plots with the spectral
measure of the finite truncated tree. To be more precise, fix $c$
and let $T_{c}$ be the infinite regular tree of degree $c$. The spectral
measure $\nu$ for the Laplacian of this tree is given by
equation (\ref{Mckay}). Consider now, for each finite $n$ the spectral
measure $\nu_{n}$ of the truncated finite tree constructed from
$T_{c}$ by just keeping only the first $n$ generations. It is a
well known result (see ~\cite{strang}) that these measures do
converge to a measure $\nu_{\infty}$. However, $\nu_{\infty}$ and
$\nu$ are very different. For instance, the measure $\nu_{\infty}$
has atoms while $\nu$ does not. This is due to the fact that repeated
eigenvalues occur with large multiplicities. The main heuristic
reason for this phenomena is that the truncated tree has a large
number of nodes with degree one creating a significant boundary
effect. See figure \ref{fig:mckay2} to see
this phenomena.

\vspace{0.3cm}
\noindent These figures clearly show that the distribution of the spectrum
of large Erd\"os--Renyi random graphs provide a better approximation for the
spectral measure of the corresponding infinite regular tree in the bulk
region than do large finite truncated trees of the same degree.  We do not
yet have a complete explanation for this.  We note that this result is
in contrast to the regime $n p_n \to \infty$ where the distribution 
of the eigenvalues follow the well-known semi-circle law \cite{fan2}.  

\vspace{0.3cm}
\noindent Of course the two spectral measures $\mu_{\infty}$ and $\nu_\infty$ are not exactly
the same; we have shown already that $\nu_\infty$ does not have a spectral
gap, whereas $\mu_{\infty}$ does. We also observe in figure~\ref{fig:mckay1} that the measure $\mu_{n}$
for small values of $c$ seems to have a spike at 1.  As seen in
figure~\ref{fig:mckay2}, the size of this spike seems to decrease
as $c$ increases, but we do not know if the spike disappears as
$n\to\infty$ for {\it fixed\/} small $c$. Nevertheless, the close
similarity observed between $\mu_{\infty}$ and $\nu_\infty$ naturally raises
the question: what is the probability distribution of the measure
$\mu_\infty$? More generally, if we consider the branching process
generated by any probability distribution in the natural numbers
$\N$, what is the spectral measure of the normalized Laplacian for
this graph?

\vspace{0.5cm}
\noindent {\it Acknowledgement}. This work was supported by AFOSR Grant No. FA9550-08-1-0064.

\end{document}